\documentclass[twoside,twocolumn,9pt]{article}
\usepackage{extsizes}
\usepackage[super,sort&compress,comma]{natbib} 
\usepackage[version=3]{mhchem}
\usepackage[left=1.5cm, right=1.5cm, top=1.785cm, bottom=2.0cm]{geometry}
\usepackage{balance}
\usepackage{mathptmx}
\usepackage{sectsty}
\usepackage{graphicx} 
\usepackage{lastpage}
\usepackage[format=plain,justification=justified,singlelinecheck=false,font={stretch=1.125,small,sf},labelfont=bf,labelsep=space]{caption}
\usepackage{float}
\usepackage{fancyhdr}
\usepackage{fnpos}
\usepackage[english]{babel}
\addto{\captionsenglish}{%
  
}
\usepackage{commath,amsmath,amssymb}
\usepackage{array}
\usepackage{droidsans}
\usepackage{charter}
\usepackage[T1]{fontenc}
\usepackage[usenames,dvipsnames]{xcolor}
\usepackage{setspace}
\usepackage[compact]{titlesec}
\usepackage{hyperref}

\usepackage{epstopdf}

\definecolor{cream}{RGB}{222,217,201}
\title{Analysis of Disordered Trusses using Network Laplacian}
\author{\large{Sean Fancher,$^{\ddag}$\textit{$^{a,\P}$} Niranjan Sarpangala,\textit{$^{b,\P}$} Prashant K. Purohit,\textit{$^{c}$} and Eleni Katifori\textit{$^{b,d}$}} }
\date{}

\begin{document}

\pagestyle{fancy}
\thispagestyle{plain}
\fancypagestyle{plain}{
\renewcommand{\headrulewidth}{0pt}
}

\makeFNbottom
\makeatletter
\renewcommand\LARGE{\@setfontsize\LARGE{15pt}{17}}
\renewcommand\Large{\@setfontsize\Large{12pt}{14}}
\renewcommand\large{\@setfontsize\large{10pt}{12}}
\renewcommand\footnotesize{\@setfontsize\footnotesize{7pt}{10}}
\makeatother

\renewcommand{\thefootnote}{\fnsymbol{footnote}}
\renewcommand\footnoterule{\vspace*{1pt}%
\color{cream}\hrule width 3.5in height 0.4pt \color{black}\vspace*{5pt}} 
\setcounter{secnumdepth}{5}

\makeatletter 
\renewcommand\@biblabel[1]{#1}            
\renewcommand\@makefntext[1]%
{\noindent\makebox[0pt][r]{\@thefnmark\,}#1}
\makeatother 
\renewcommand{\figurename}{\small{Fig.}~}
\sectionfont{\sffamily\Large}
\subsectionfont{\normalsize}
\subsubsectionfont{\bf}
\setstretch{1.125} 
\setlength{\skip\footins}{0.8cm}
\setlength{\footnotesep}{0.25cm}
\setlength{\jot}{10pt}
\titlespacing*{\section}{0pt}{4pt}{4pt}
\titlespacing*{\subsection}{0pt}{15pt}{1pt}

\fancyfoot{}
\fancyfoot[RO]{\thepage}
\fancyfoot[LE]{\thepage}
\fancyhead{}
\renewcommand{\headrulewidth}{0pt} 
\renewcommand{\footrulewidth}{0pt}
\setlength{\arrayrulewidth}{1pt}
\setlength{\columnsep}{6.5mm}
\setlength\bibsep{1pt}

\makeatletter 
\newlength{\figrulesep} 
\setlength{\figrulesep}{0.5\textfloatsep} 

\newcommand{\topfigrule}{\vspace*{-1pt}%
\noindent{\color{cream}\rule[-\figrulesep]{\columnwidth}{1.5pt}} }

\newcommand{\botfigrule}{\vspace*{-2pt}%
\noindent{\color{cream}\rule[\figrulesep]{\columnwidth}{1.5pt}} }

\newcommand{\dblfigrule}{\vspace*{-1pt}%
\noindent{\color{cream}\rule[-\figrulesep]{\textwidth}{1.5pt}} }

\makeatother


\maketitle

\section*{Abstract} 
\textbf{Truss structures, with distributed mass elements, at macro-scale are common in a number of engineering applications and are now being increasingly used at the micro-scale to construct metamaterials. In analyzing the properties of a given truss structure, it is often necessary to understand how stress waves propagate through the system and/or its dynamic modes under time dependent loading so as to allow for maximally efficient use of space and material. This can be a computationally challenging task for particularly large or complex structures, with current methods requiring fine spatial discretization or evaluations of sizable matrices. Here we present a spectral method to compute the dynamics of trusses inspired by results from fluid flow networks. Our model accounts for the full dynamics of linearly elastic truss elements via a network Laplacian; a matrix object which couples the motions of the structure joints. We show that this method is equivalent to the continuum limit of linear finite element methods as well as capable of reproducing natural frequencies and modes determined by more complex and computationally costlier methods. Our results show that balls-and-springs models inadequately describe dynamics, especially at short times relative to wave propagation time through rods. Furthermore, we illustrate the method's utility in optimizing target joint displacements using impedance matching and resonance-based schemes, offering a computationally efficient approach for analyzing large, complex truss structures.}\\


\renewcommand*\rmdefault{bch}\normalfont\upshape
\rmfamily
\section*{}
\vspace{-1cm}


\footnotetext{\textit{a}~Department of Biophysics, University of Michigan, Ann Arbor, MI}
\footnotetext{\textit{b}~Department of Physics and Astronomy, University of Pennsylvania, Philadelphia, PA}
\footnotetext{\textit{c}~Department of Mechanical Engineering and Applied Mechanics, University of Pennsylvania, Philadelphia PA}
\footnotetext{\textit{d}~Center for Computational Biology, Flatiron Institute, New York, NY}

\footnotetext{\P~These authors contributed equally to the work}


\section{Introduction}
\label{sec:intro}
Trusses have been a mainstay of structural engineering for a substantial amount of human history, with applications including bridges, buildings, airplanes, and spacecraft. While these applications at the scale of several tens of meters are well established, truss metamaterials with microstructures in the range of a few millimeters are of intense research interest currently due to their easy manufacturability by 3D printing and other techniques \cite{askari2020additive}. They are also relevant in understanding the properties and adaptation of biological systems ranging from cytoskeleton to bone tissues \cite{torres2019bone, van2020mechanoresponse} in animal skeletons \cite{jung2018comparative}. However, efficiently analyzing the dynamic behavior of complex truss structures, especially in disordered or hierarchical configurations, remains a significant challenge in the field of soft matter physics.\\

Analysing wave propagation through truss metamaterials \cite{mueller2019energy,glaesener2021viscoelastic}, 
is particularly interesting since they can be designed to allow only certain wavelengths/frequencies to propagate, thus enabling applications in acoustic absorbers and transmitters \cite{wang2015locally,krushynska2018accordion}. Similar techniques have also been applied to structures comprised of polymer networks for the purpose of understanding response to propagating loads and defect particles \cite{zhang2015novel,lu2022double}. These applications involve precise determination of the dynamic behavior of the structure. Furthermore, current imaging techniques for the measurement of local strains have become so sophisticated that it is possible to observe wavefronts propagating through individual elements at the microscale \cite{branch2017controlling}. These experimental developments enable the validation of detailed computational models which was not possible before.\\

The dynamic behavior of truss structures is most commonly studied using the finite element method. Each member of the structure is discretized into appropriate truss elements, then stiffness and mass matrices are assembled, then the equations of motion (including external forces) are written in matrix form starting from Newton's laws. The natural frequencies and mode shapes (which are critical in structural design and in understanding the dispersion relations of metamaterials) are computed by solving an eigenvalue problem using the stiffness and mass matrices \cite{cook2007concepts}. The frequencies obtained depend on the mass and stiffness matrices; they vary depending on how fine the discretization is, whether the consistent or lumped mass matrix (or a combination) is used, and on the shape functions used to assemble the stiffness and mass matrices. Often, very fine discretization is required to get an accurate estimate of the natural frequencies and the computational costs can be prohibitive if the structure has a large number of degrees of freedom. Methods that retain the full element dynamics have also been developed via solutions in Fourier space (reverberation matrix method) \cite{howard1998analysis, pao1999dynamic} or integral operators \cite{polz2019wave}. These produce accurate solutions for large structures but require more sophisticated mathematics and analysis algorithms. While it is possible to directly simulate the temporal response of a truss structure by monitoring the propagation of reactive wavefronts \cite{howard1998analysis,messner2015wave,trainiti2016wave}, this can become computationally costly as more wavefronts are produced via scattering. \\

In this article, we present a new method for determining the dynamic behavior of a given truss structure based on our work on fluid flow networks \cite{fancher2022mechanical}. We retain a full description of the element dynamics much like the methodology of \cite{howard1998analysis,pao1999dynamic}, but we maintain a focus on the motions of the structure joints rather than stress within the elements. By assuming linear stress and displacement dynamics, we can express the propagation of stress waves in terms of the joint motions in Fourier space and develop a linear relation between the displacements of and forces acting on the joints. This allows our resulting matrix objects to be smaller in size and thus more computationally efficient to analyze in comparison to those produced via element stress calculations. This method is different from the spectral element methods \cite{jf1997wave, chakraborty2003spectrally, lee2009spectral, an2020elastic, zuo2016numerical} which are similar to finite element methods except that the shape functions are chosen to be orthogonal functions as in a Fourier basis. Thus, if the bars are discretized using spectral elements then the size of the problem scales with the number of bars. In the technique presented here the size of the problem scales with the number of joints.\\

We demonstrate the method's accuracy and efficiency through analysis of a simple structure (square frame with crossbar) and extend its application to disordered networks, revealing the inadequacy of ball-and-spring models at high frequencies. In contrast to the balls and spring method, our spectral method gives consistent results across frequency ranges. Finally, we show that our method can be used to optimize target joint motion in response to input joint oscillations, illustrating its usefulness in designing dynamics of complex, disordered systems.




\section{Methods}
\label{sec:methods}

\begin{figure}
    \includegraphics[width=0.38\textwidth]{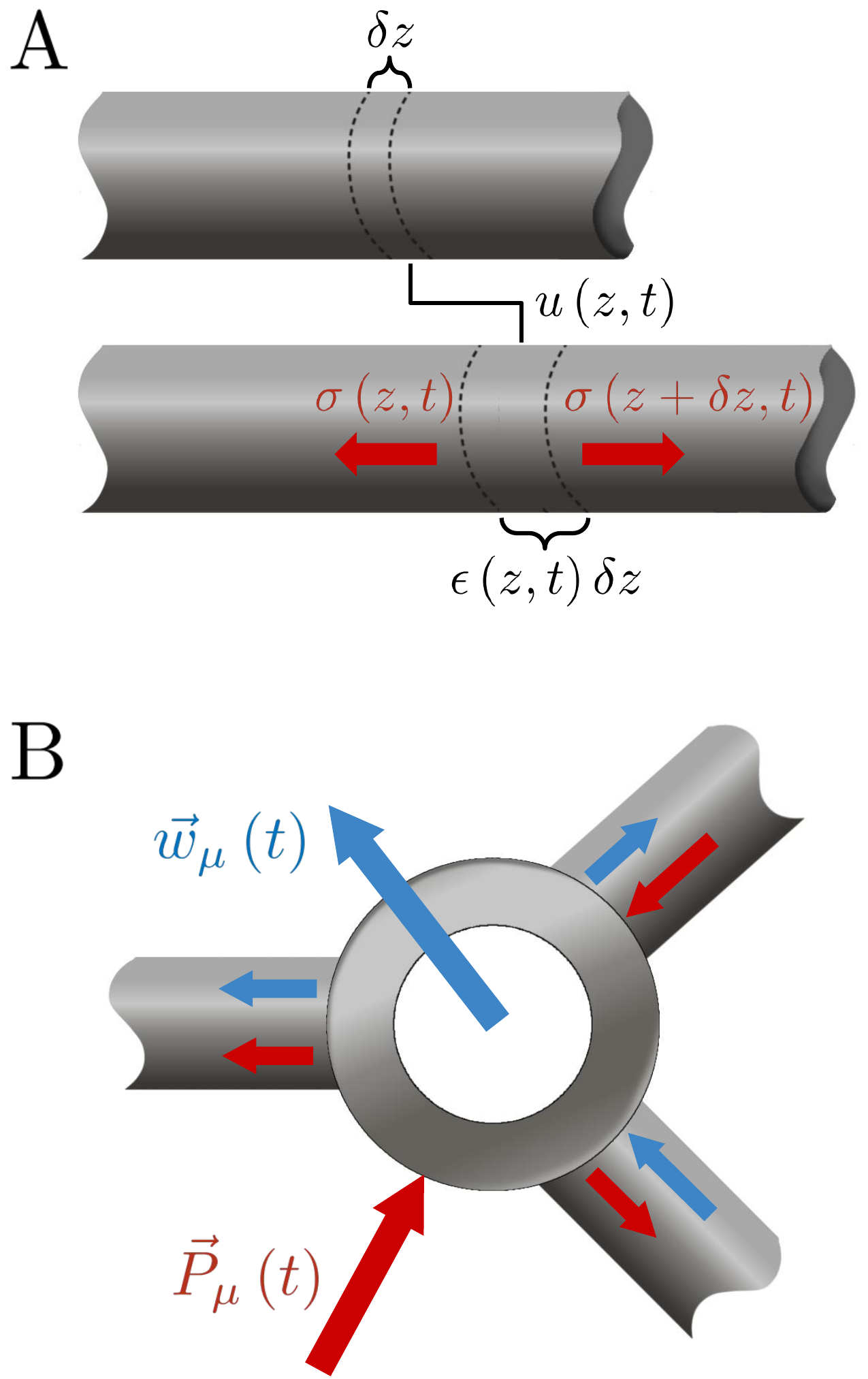}
  \caption{A) Single rods undergo position dependent displacement, $u(z,t)$, which causes strain, $\epsilon(z,t)$, within the material. This results in a buildup of stress, $\sigma(z,t)$, which resists the strain. B) Multiple rods can be connected into a network via joints. Each joint can move with velocity $\vec{w}_{\mu}(t)$ and induce similar velocity within the individual elements (blue arrows). Similarly, a force, $\vec{P}_{\mu}(t)$, can be applied to the joint but will be resisted by the elements (red arrows).}
  \label{Fig:cartoons}
\end{figure}

\subsection{Single Rod Dynamics}
\label{sec:SRD}

We begin by considering a solid, cylindrical rod of length $L$ and cross sectional area $A$. As the rod undergoes tension and compression in response to externally applied forces, each infinitesimally thin segment is displaced from its rest position by an amount $u(z,t)$, where $z$ is the distance from the $z=0$ end of the rod in its unperturbed configuration and $t$ is time. Specifically, $u(z,t)>0$ implies the segment has shifted in the direction of increasing $z$ and $u(z,t)<0$ implies a shift in the direction of decreasing $z$, as depicted in Fig. \ref{Fig:cartoons} (A).

From this displacement field we can derive two other critically important fields; the velocity field, $v(z,t)$, and strain field, $\epsilon (z,t)$. Defining the sign of $v$ and $\epsilon$ to denote movement in the direction of $z$ and tensile expansion of the material respectively allows for the relations

\begin{equation}
v\left(z,t\right) = \frac{\partial u}{\partial t}, \quad\quad\quad \epsilon\left(z,t\right) = \frac{\partial u}{\partial z}.
\label{vepdef}
\end{equation}

\noindent Under the assumption that the material is linearly elastic, the stress field, $\sigma(z,t)$, can be expressed simply as $\sigma(z,t)=E\epsilon(z,t)$, where $E$ is the Young's modulus. This forces the sign of $\sigma$ to be such that $\sigma>0$ represents tensile stress while $\sigma<0$ represents compressive stress. These definitions of $v$, $\epsilon$, and $\sigma$ immediately allow for the identity

\begin{equation}
\frac{\partial v}{\partial z}-\frac{1}{E}\frac{\partial\sigma}{\partial t} = \frac{\partial}{\partial z}\left(\frac{\partial u}{\partial t}\right)-\frac{\partial}{\partial t}\left(\frac{\partial u}{\partial z}\right) = 0.
\label{masscon}
\end{equation}

We now consider the forces acting on an infinitesimally thin segment of the rod, such as that of width $\delta z$ shown in Fig. \ref{Fig:cartoons} (A). Following our sign convention for the stress field, Newton's second law takes the form

\begin{equation}
\frac{\partial\sigma}{\partial z}-\rho\frac{\partial v}{\partial t} = 0.
\label{forcebalance}
\end{equation}

\noindent where $\rho$ is the material density and the $\delta z\to 0$ has been taken. Together, Eqs. \ref{masscon} and \ref{forcebalance} represent the dynamic coupling between the stress and velocity fields of the rod \cite{howard1998analysis,pao1999dynamic}. Additionally, we can express these in terms of the displacement field to transform Eq. \ref{forcebalance} into

\begin{equation}
E\frac{\partial^{2}u}{\partial z^{2}}-\rho\frac{\partial^{2}u}{\partial t^{2}} = 0,
\label{dispwave}
\end{equation}

\noindent thus producing the simple wave equation wherein waves can propagate through the displacement field at speed $c=\sqrt{E/\rho}$.

One important implication of Eq. \ref{dispwave} is that $u(z,t)$ can be expanded into its Fourier modes via

\begin{equation}
u\left(z,t\right) = \int d\omega\>\left(F\left(\omega\right)e^{i\omega\left(t-\frac{z}{c}\right)}+B\left(\omega\right)e^{i\omega\left(t+\frac{z}{c}\right)}\right),
\label{uFint}
\end{equation}

\noindent where $F(\omega)$ and $B(\omega)$ are the forward and backward wave amplitudes respectively and the integral being over all real $\omega$ is implied. This in turn allows for the velocity and stress fields to be expressed as

\begin{subequations}
\begin{equation}
v\left(z,t\right) = \frac{\partial u}{\partial t} = \int d\omega\>i\omega\left(F\left(\omega\right)e^{i\omega\left(t-\frac{z}{c}\right)}+B\left(\omega\right)e^{i\omega\left(t+\frac{z}{c}\right)}\right),
\label{vFint}
\end{equation}

\begin{equation}
\sigma\left(z,t\right) = E\frac{\partial u}{\partial z} = \Gamma\int d\omega\>i\omega\left(-F\left(\omega\right)e^{i\omega\left(t-\frac{z}{c}\right)}+B\left(\omega\right)e^{i\omega\left(t+\frac{z}{c}\right)}\right),
\label{sFint}
\end{equation}
\label{vsFint}
\end{subequations}

\noindent where $\Gamma=E/c=\sqrt{E\rho}$ is the material impedance.

\subsection{Rod and Joint Network Construction}
\label{sec:RJNC}

We now consider a network comprised of rods obeying the dynamic equations outlined thus far connected via a series of joints such as the one depicted in Fig. \ref{Fig:cartoons} (B). Here we will use the index $\mu$ to denote a particular joint and the index $\nu$ to denote a joint connected to $\mu$ through one of the network rods. In this way, $\nu\in\mathcal{N}_{\mu}$, where $\mathcal{N}_{\mu}$ is the set of all joints connected to joint $\mu$ through a single rod. The rods themselves and their various fields will be labelled with a two component index, $\mu\nu$, comprised of the two joints the rod connects. Specifically, $u_{\mu\nu}(z,t)$ is the displacement field of rod $\mu\nu$ with the $z=0$ end of the rod being at joint $\mu$. Similar notation also applies to other fields as well as rod specific parameters such as the length, $L_{\mu\nu}$, and cross sectional area, $A_{\mu\nu}$. Exchanging the index order thus also reverses the directionality of the rod, leading to the relations

\begin{subequations}
\begin{equation}
u_{\mu\nu}\left(z,t\right) = -u_{\nu\mu}\left(L_{\mu\nu}-z,t\right),
\label{uexchange}
\end{equation}
\begin{equation}
v_{\mu\nu}\left(z,t\right) = -v_{\nu\mu}\left(L_{\mu\nu}-z,t\right),
\label{vexchange}
\end{equation}
\begin{equation}
\sigma_{\mu\nu}\left(z,t\right) = \sigma_{\nu\mu}\left(L_{\mu\nu}-z,t\right).
\label{sexchange}
\end{equation}
\label{uvsexchange}
\end{subequations}

\noindent We can also apply the Fourier expansion used in Eq. \ref{uFint} alongside this notation to obtain the index exchange laws for $F_{\mu\nu}(\omega)$ and $B_{\mu\nu}(\omega)$. Letting $\tau_{\mu\nu}=L_{\mu\nu}/c_{\mu\nu}$ allows for these to be expressed as

\begin{equation}
F_{\mu\nu}\left(\omega\right) = -B_{\nu\mu}\left(\omega\right)e^{i\omega\tau_{\mu\nu}}, \quad\quad\quad B_{\mu\nu}\left(\omega\right) = -F_{\nu\mu}\left(\omega\right)e^{-i\omega\tau_{\mu\nu}}.
\label{FBexchange}
\end{equation}

\noindent Finally, we can Fourier transform Eq. \ref{vsFint} in time and combine the result with Eq. \ref{FBexchange} to produce the relations

\begin{subequations}
\begin{equation}
\tilde{v}_{\mu\nu}\left(0,\omega\right) = i\omega\left(F_{\mu\nu}\left(\omega\right)+B_{\mu\nu}\left(\omega\right)\right),
\label{v0FT}
\end{equation}
\begin{equation}
\tilde{\sigma}_{\mu\nu}\left(0,\omega\right) = i\omega\Gamma_{\mu\nu}\left(-F_{\mu\nu}\left(\omega\right)+B_{\mu\nu}\left(\omega\right)\right),
\label{s0FT}
\end{equation}
\label{vsFT}
\end{subequations}

\begin{subequations}
\begin{align}
\tilde{v}_{\mu\nu}\left(0,\omega\right) = -\tilde{v}_{\nu\mu}\left(0,\omega\right)\cos\left(\omega\tau_{\mu\nu}\right)-i\Gamma_{\mu\nu}^{-1}\tilde{\sigma}_{\nu\mu}\left(0,\omega\right)\sin\left(\omega\tau_{\mu\nu}\right),
\label{vFTexchange}
\end{align}
\begin{equation}
\tilde{\sigma}_{\mu\nu}\left(0,\omega\right) =  \tilde{\sigma}_{\nu\mu}\left(0,\omega\right)\cos\left(\omega\tau_{\mu\nu}\right)+i\Gamma_{\mu\nu}\tilde{v}_{\nu\mu}\left(0,\omega\right)\sin\left(\omega\tau_{\mu\nu}\right).
\label{sFTexchange}
\end{equation}
\label{vsFTexchange}
\end{subequations}

From here we assign a $D_{\mu}$-dimensional coordinate system, denoted as $\mathbb{R}_{\mu}^{D}$, to the $\mu$th joint such that the origin is located at the rest location of the joint and the set of rods connected to that joint span $\mathbb{R}_{\mu}^{D}$. For example, if two rods are connected at 180° (antiparallel), the dimension of the joint is 1.  On the other hand, if they are not parallel or antiparallel, say the rods are connected at an angle $60^o$, then the dimension of the joint is 2. We can then define a unit vector, $\hat{e}_{\mu\nu}\in\mathbb{R}_{\mu}^{D}$, to rod $\mu\nu$ with equivalent directionality, thus implying that $\hat{e}_{\mu\nu}$ points from joint $\mu$ to joint $\nu$. Of note is that since the opposing vector $\hat{e}_{\nu\mu}$ exists in $\mathbb{R}_{\nu}^{D}$, there is no implicit index exchange relation between $\hat{e}_{\mu\nu}$ and $\hat{e}_{\nu\mu}$ without first defining the relation between $\mathbb{R}_{\mu}^{D}$ and $\mathbb{R}_{\nu}^{D}$. However, if $\hat{e}_{\mu\nu}$ and $\hat{e}_{\nu\mu}$ are expressed in the global coordinate system, denoted $\mathbb{R}_{g}^{D}$, then they must of course point in opposing directions. This is particularly easy to achieve if $\text{dim}(\mathbb{R}_{\mu}^{D})=\text{dim}(\mathbb{R}_{g}^{D})$ for all $\mu$. In this case, we can define an ``aligned coordinate set" in which $\vec{x}_{\mu}=\vec{x}_{g}-\vec{r}_{\mu}$; where $\vec{x}_{\mu}$ is a position vector in $\mathbb{R}_{\mu}^{D}$, $\vec{x}_{g}$ is the same position in $\mathbb{R}_{g}^{D}$, and $\vec{r}_{\mu}$ is the position of the $\mu$th joint in $\mathbb{R}_{g}^{D}$.

With the coordinate systems defined, we next investigate the dynamics of the joints by assuming that each joint is massless and incapable of carrying force. Newton's second law applied to joint $\mu$ then takes the form

\begin{equation}
\vec{P}_{\mu}\left(t\right)+\sum_{\nu\in\mathcal{N}_{\mu}}\hat{e}_{\mu\nu}A_{\mu\nu}\sigma_{\mu\nu}\left(0,t\right) = \vec{0},
\label{jointforce}
\end{equation}

\noindent where $\vec{P}_{\mu}\left(t\right)$ is the force being applied to the joint by some entity external to the system, again expressed within $\mathbb{R}_{\mu}^{D}$. Finally, the joint itself moves within this coordinate system with a velocity given by the vector $\vec{w}_{\mu}(t)$. Enforcing that the $z=0$ end of rod $\mu\nu$ must have a velocity equivalent to the projection of $\vec{w}_{\mu}$ in the direction of $\hat{e}_{\mu\nu}$ yields the condition

\begin{equation}
\vec{w}_{\mu}\cdot\hat{e}_{\mu\nu} = v_{\mu\nu}\left(0,t\right) \quad\forall\quad \nu\in\mathcal{N}_{\mu}.
\label{jointvelcon}
\end{equation}

Eqs. \ref{jointforce} and \ref{jointvelcon} provide the connectivity laws of the network and define how the stress and velocity fields of different rods interact \cite{howard1998analysis,pao1999dynamic}. Their form can be somewhat simplified by introducing the matrix $\overset\Leftrightarrow{e}_{\mu}$, defined to be a $D_{\mu}\times|\mathcal{N}_{\mu}|$ matrix in which each column is a distinct $\hat{e}_{\mu\nu}$. The joint stress and velocity vectors, $\vec{\sigma}_{\mu}(t)$ and $\vec{v}_{\mu}(t)$, can then be defined as column vectors of length $|\mathcal{N}_{\mu}|$ such that their respective $j$th components give $\sigma_{\mu\nu}(0,t)$ and $v_{\mu\nu}(0,t)$ for the same $\nu$ as was used to generate the $j$th column of $\overset\Leftrightarrow{e}_{\mu}$. Similarly, we construct $\overset\Leftrightarrow{A}_{\mu}$ as a diagonal matrix of size $|\mathcal{N}_{\mu}|\times|\mathcal{N}_{\mu}|$ whose $j$th diagonal entry is $A_{\mu\nu}$. Finally, treating $\vec{P}_{\mu}\left(t\right)$ and $\vec{w}_{\mu}$ as a column vectors allows Eqs. \ref{jointforce} and \ref{jointvelcon} to be expressed as

\begin{subequations}
\begin{equation}
\overset\Leftrightarrow{e}_{\mu}\overset\Leftrightarrow{A}_{\mu}\vec{\sigma}_{\mu} = -\vec{P}_{\mu},
\label{jointforcemat}
\end{equation}
\begin{equation}
\left(\vec{w}_{\mu}\right)^{\text{T}}\overset\Leftrightarrow{e}_{\mu} = \left(\vec{v}_{\mu}\right)^{\text{T}}.
\label{jointvelconmat}
\end{equation}
\label{jointmateqs}
\end{subequations}

Of note is that due to our construction of $\mathbb{R}_{\mu}^{D}$ we have $D_{\mu}\le |\mathcal{N}_{\mu}|$. In the specific case of equality, $\overset\Leftrightarrow{e}_{\mu}$ must be an invertible square matrix due to the set of vectors $\{\hat{e}_{\mu\nu}\> |\>\nu\in\mathcal{N}_{\mu}\}$ spanning the local coordinate system. This causes Eqs. \ref{jointforcemat} and \ref{jointvelconmat} to become a bijective linear relation between the rod dynamics, $\vec{\sigma}_{\mu}$ and $\vec{v}_{\mu}$, and the joint dynamics, $\vec{P}_{\mu}$ and $\vec{w}_{\mu}$. Thus, the rods are all effectively uncoupled from each other when $D_{\mu}=|\mathcal{N}_{\mu}|$. This is because for any given stress and velocity field within a single rod there must exist an applied force and joint velocity such that Eqs. \ref{jointforcemat} and \ref{jointvelconmat} are satisfied without inducing any additional dynamics in any other rod.

\subsection{The Laplacian Block Matrix}
\label{sec:MML}

We now seek to construct a global relation between the motions of every joint in the network. We achieve this by first noting the similarities between the element velocity and stress fields with the flow velocity and pressure of a fluid in a compliant vessel, then following a methodology previously developed for calculating the pressure distribution in such a fluid flow network \cite{fancher2022mechanical}. To begin, we combine Eqs. \ref{vexchange} and \ref{v0FT} to produce

\begin{subequations}
\begin{align}
\tilde{v}_{\nu\mu}\left(0,\omega\right) & = -\tilde{v}_{\mu\nu}\left(L_{\mu\nu},\omega\right) \nonumber \\ & = -i\omega\left(F_{\mu\nu}\left(\omega\right)e^{-i\omega\tau_{\mu\nu}}+B_{\mu\nu}\left(\omega\right)e^{i\omega\tau_{\mu\nu}}\right)
\end{align}
\label{vLtoFB}
\end{subequations}

\noindent Solving $F_{\mu\nu}$ and $B_{\mu\nu}$ yields

\begin{subequations}
\begin{align}
F_{\mu\nu}\left(\omega\right) &= \frac{\tilde{v}_{\mu\nu}\left(0,\omega\right)+\tilde{v}_{\nu\mu}\left(0,\omega\right)e^{-i\omega\tau_{\mu\nu}}}{i\omega\left(1-e^{-2i\omega\tau_{\mu\nu}}\right)} \nonumber \\
&= \frac{\tilde{\vec{w}}_{\mu}^{\text{T}}\left(\omega\right)\hat{e}_{\mu\nu}+\tilde{\vec{w}}_{\nu}^{\text{T}}\left(\omega\right)\hat{e}_{\nu\mu}e^{-i\omega\tau_{\mu\nu}}}{i\omega\left(1-e^{-2i\omega\tau_{\mu\nu}}\right)}
\label{Ftov0L}
\end{align}
\begin{align}
B_{\mu\nu}\left(\omega\right) &= \frac{\tilde{v}_{\mu\nu}\left(0,\omega\right)+\tilde{v}_{\nu\mu}\left(0,\omega\right)e^{i\omega\tau_{\mu\nu}}}{i\omega\left(1-e^{2i\omega\tau_{\mu\nu}}\right)} \nonumber \\&= \frac{\tilde{\vec{w}}_{\mu}^{\text{T}}\left(\omega\right)\hat{e}_{\mu\nu}+\tilde{\vec{w}}_{\nu}^{\text{T}}\left(\omega\right)\hat{e}_{\nu\mu}e^{i\omega\tau_{\mu\nu}}}{i\omega\left(1-e^{2i\omega\tau_{\mu\nu}}\right)}
\label{Btov0L}
\end{align}
\label{FBtov0L}
\end{subequations}

\noindent where we have replaced $\tilde{v}_{\mu\nu}$ and $\tilde{v}_{\nu\mu}$ with the Fourier transformed velocities of joints $\mu$ and $\nu$ respectively as per Eq. \ref{jointvelcon}. From here we use Eq. \ref{s0FT} to express the Fourier transformed stress as

\begin{align}
\tilde{\sigma}_{\mu\nu}\left(0,\omega\right) &= i\omega\Gamma_{\mu\nu}\left(-F_{\mu\nu}\left(\omega\right)+B_{\mu\nu}\left(\omega\right)\right) \nonumber \\
&= \Gamma_{\mu\nu}\left(-\frac{\tilde{\vec{w}}_{\mu}^{\text{T}}\left(\omega\right)\hat{e}_{\mu\nu}+\tilde{\vec{w}}_{\nu}^{\text{T}}\left(\omega\right)\hat{e}_{\nu\mu}e^{-i\omega\tau_{\mu\nu}}}{1-e^{-2i\omega\tau_{\mu\nu}}}\right.\\& \left. \qquad +\frac{\tilde{\vec{w}}_{\mu}^{\text{T}}\left(\omega\right)\hat{e}_{\mu\nu}+\tilde{\vec{w}}_{\nu}^{\text{T}}\left(\omega\right)\hat{e}_{\nu\mu}e^{i\omega\tau_{\mu\nu}}}{1-e^{2i\omega\tau_{\mu\nu}}}\right) \nonumber\\
&= \frac{i\Gamma_{\mu\nu}}{\sin\left(\omega\tau_{\mu\nu}\right)}\left(\hat{e}_{\mu\nu}^{\text{T}}\tilde{\vec{w}}_{\mu}\left(\omega\right)\cos\left(\omega\tau_{\mu\nu}\right)+\hat{e}_{\nu\mu}^{\text{T}}\tilde{\vec{w}}_{\nu}\left(\omega\right)\right)
\label{s0towmn}
\end{align}

\noindent Inserting Eq. \ref{s0towmn} into the Fourier transform of Eq. \ref{jointforce} and introducing the parameter $\Lambda_{\mu\nu}=A_{\mu\nu}\Gamma_{\mu\nu}$ then yields

\begin{align}
-\tilde{\vec{P}}_{\mu}\left(\omega\right) &= \sum_{\nu\in\mathcal{N}_{\mu}}\hat{e}_{\mu\nu}A_{\mu\nu}\tilde{\sigma}_{\mu\nu}\left(0,\omega\right) \nonumber \\ &= \sum_{\nu\in\mathcal{N}_{\mu}}\hat{e}_{\mu\nu}\frac{i\Lambda_{\mu\nu}}{\sin\left(\omega\tau_{\mu\nu}\right)}\left(\hat{e}_{\mu\nu}^{\text{T}}\tilde{\vec{w}}_{\mu}\left(\omega\right)\cos\left(\omega\tau_{\mu\nu}\right)+\hat{e}_{\nu\mu}^{\text{T}}\tilde{\vec{w}}_{\nu}\left(\omega\right)\right)
\label{s0sumtowmn}
\end{align}

Based on Eq. \ref{s0sumtowmn} we can construct the block vectors $\mathbf{\vec{W}}$ and $\mathbf{\vec{P}}$ as well as the block matrix $\mathbf{\overset\Leftrightarrow{D}}$. These are objects whose individual components are themselves vectors and matrices. Specifically, $\mathbf{\vec{W}}$ and $\mathbf{\vec{P}}$ have $J$ components each, where $J$ is the number of joints in the network, with the $\mu$th components being the vectors $\tilde{\vec{w}}_{\mu}$ and $\tilde{\vec{P}}_{\mu}$ respectively, each expressed in terms of $\mathbb{R}_{\mu}^{D}$. Similarly, $\mathbf{\overset\Leftrightarrow{D}}$, denoted here as the network Laplacian, has a $J\times J$ structure with components defined as

\begin{equation}
\mathbf{\overset\Leftrightarrow{D}}_{\mu\nu} = \begin{cases}
\sum_{\gamma\in\mathcal{N}_{\mu}}\Lambda_{\mu\gamma}\omega\cot\left(\omega\tau_{\mu\gamma}\right)\hat{e}_{\mu\gamma}\hat{e}_{\mu\gamma}^{\text{T}} & \mu=\nu \\
\Lambda_{\mu\nu}\omega\csc\left(\omega\tau_{\mu\nu}\right)\hat{e}_{\mu\nu}\hat{e}_{\nu\mu}^{\text{T}} & \nu\in\mathcal{N}_{\mu} \\
\overset\Leftrightarrow{0} & \text{otherwise} \end{cases},
\label{Ddef}
\end{equation}

\noindent thus making the $\mu\nu$ component of $\mathbf{\overset\Leftrightarrow{D}}$ a $D_{\mu}\times D_{\nu}$ matrix that transforms a vector in $\mathbb{R}_{\nu}^{D}$ into one in $\mathbb{R}_{\mu}^{D}$. Given these definitions, Eq. \ref{s0sumtowmn} clearly dictates

\begin{equation}
-\frac{1}{i\omega}\mathbf{\overset\Leftrightarrow{D}}\mathbf{\vec{W}} = -\mathbf{\vec{P}} \quad\quad\quad \implies \quad\quad\quad \mathbf{\overset\Leftrightarrow{D}}\mathbf{\vec{U}} = \mathbf{\vec{P}}.
\label{DWPrel}
\end{equation}

Eq. \ref{DWPrel} provides a direct relation between the force applied to the system and its dynamics. The first equality generates a linear transformation between the joint velocities and forces. The second equality generates a similar relation in terms of $\mathbf{\vec{U}}$, the block vector of Fourier transformed joint displacements, and is obtained from the first by noting that $\mathbf{\vec{W}}=i\omega\mathbf{\vec{U}}$ since velocity is the time derivative of displacement. Of note is that the network Laplacian derived here is similar in function to the global scattering matrix of the system \cite{pao1999dynamic}, but importantly represents a relation over the joints rather than elements. Thus, the network Laplacian is typically smaller in size and more computationally manageable.

\section{Results}
\subsection{Stiffness and mass matrix method is just a second order approximation of our truss network method}
\label{sec:SMM}

Our theory, in particular Eq. \ref{DWPrel} can be interpreted as a dynamic extension of the static stiffness matrix method. In this method, $\mathbf{\overset\Leftrightarrow{K}}$ is the static stiffness matrix, $\mathbf{\vec{U}}$ is the displacement vector whose elements give the static displacement from equilibrium of each joint, and $\mathbf{\vec{Q}}$ is the vector of externally applied forces that maintain this out of equilibrium positioning. The stiffness matrix can be defined by treating each rod as a Hookian spring of stiffness $k_{\mu\nu}=A_{\mu\nu}E_{\mu\nu}/L_{\mu\nu}=\Lambda_{\mu\nu}/\tau_{\mu\nu}$ so that the force it exerts at joint $\mu$ is given by $-k_{\mu\nu}\hat{e}_{\mu\nu}(\hat{e}_{\mu\nu}^{\text{T}}\vec{u}_{\mu}-\hat{e}_{\nu\mu}^{\text{T}}\vec{u}_{\nu})$, where $\vec{u}_{\mu}\in\mathbb{R}_{\mu}^{D}$ is the displacement of the $\mu$th joint in its own local coordinate system. Summing this effect over all rods and using our construction of $\mathbf{\overset\Leftrightarrow{D}}$ given by Eq. \ref{Ddef} automatically provides the relation

\begin{equation}
 \mathbf{\overset\Leftrightarrow{K}} = \lim_{\omega\to 0}\mathbf{\overset\Leftrightarrow{D}}\left(\omega\right).
\label{KDrel}
\end{equation}

\noindent From here, we can use the fact that $\mathbf{\vec{U}}$ and $\mathbf{\vec{Q}}$ represent static quantities to express their Fourier transforms as simply the vectors themselves multiplied by a $\delta$-function. This allows the second form of Eq. \ref{DWPrel} to be expressed as

\begin{equation}
\mathbf{\vec{Q}}\delta\left(\omega\right) = \mathbf{\vec{P}} = \mathbf{\overset\Leftrightarrow{D}}\mathbf{\vec{U}}\delta\left(\omega\right) = \left(\mathbf{\overset\Leftrightarrow{K}}\mathbf{\vec{U}}\right)\delta\left(\omega\right).
\label{Pstatic}
\end{equation}

\noindent Equating the prefactors of the $\delta$-functions on either side of Eq. \ref{Pstatic} yields the condition $\mathbf{\vec{Q}}=\mathbf{\overset\Leftrightarrow{K}}\mathbf{\vec{U}}$, which is precisely a statement of static equilibrium written in matrix form.

We can further explore the limiting behavior of Eq. \ref{DWPrel} by defining the consistent mass matrix of the system as

\begin{equation}
\mathbf{\overset\Leftrightarrow{M}}_{\mu\nu} = \begin{cases}
\sum_{\gamma\in\mathcal{N}_{\mu}}\int_{0}^{L_{\mu\gamma}}dz\>A_{\mu\gamma}\rho_{\mu\gamma}\left(N\left(\frac{z}{L_{\mu\gamma}}\right)\right)^{2}\hat{e}_{\mu\gamma}\hat{e}_{\mu\gamma}^{\text{T}} \qquad \mu=\nu \\
-\int_{0}^{L_{\mu\nu}}dz\>A_{\mu\nu}\rho_{\mu\nu}N\left(\frac{z}{L_{\mu\nu}}\right) \\
 \qquad  \qquad \left(1-N\left(\frac{z}{L_{\mu\nu}}\right)\right)\hat{e}_{\mu\nu}\hat{e}_{\nu\mu}^{\text{T}} \quad \qquad \qquad \nu\in\mathcal{N}_{\mu} \\
\overset\Leftrightarrow{0} \qquad \qquad \qquad \qquad  \qquad \qquad \qquad \qquad \quad  \text{otherwise} \end{cases},
\label{MCdef}
\end{equation}

\noindent where $N(x):\lbrack 0,1\rbrack\to\lbrack 0,1\rbrack$ is the shape function of the rod. Note the sign negation in the off-diagonal terms which takes into account the opposing directionalities of $\hat{e}_{\mu\nu}$ and $\hat{e}_{\nu\mu}$ when expressed in terms of $\mathbb{R}_{g}^{D}$. Using the standard linear shape function ($N(x)=1-x$) allows the integrals to be easily performed while also yielding the relation

\begin{equation}
\mathbf{\overset\Leftrightarrow{M}} = -\frac{1}{2}\lim_{\omega\to 0}\left(\frac{\partial^{2}}{\partial\omega^{2}}\left(\mathbf{\overset\Leftrightarrow{D}}\right)\right) \implies\mathbf{\overset\Leftrightarrow{D}} = \mathbf{\overset\Leftrightarrow{K}}-\omega^{2}\mathbf{\overset\Leftrightarrow{M}}+\mathcal{O}\left(\omega^{4}\right)
\label{MDrel}
\end{equation}

\noindent where $\mathbf{\overset\Leftrightarrow{D}}$ has been approximated by its Taylor expansion to second order. This expansion shows that the practice \cite{cook2007concepts} of finding the natural frequencies of the system by observing where $\det(\mathbf{\overset\Leftrightarrow{K}}-\omega^{2}\mathbf{\overset\Leftrightarrow{M}})=0$ is merely a second order approximation of defining the natural frequencies by where $\det(\mathbf{\overset\Leftrightarrow{D}})=0$.

\subsection{Truss network method gives accurate results at a lower computational cost}
\label{sec:SRD}

\begin{figure*}[h!]
    \centering
    \includegraphics[width=0.8\linewidth]{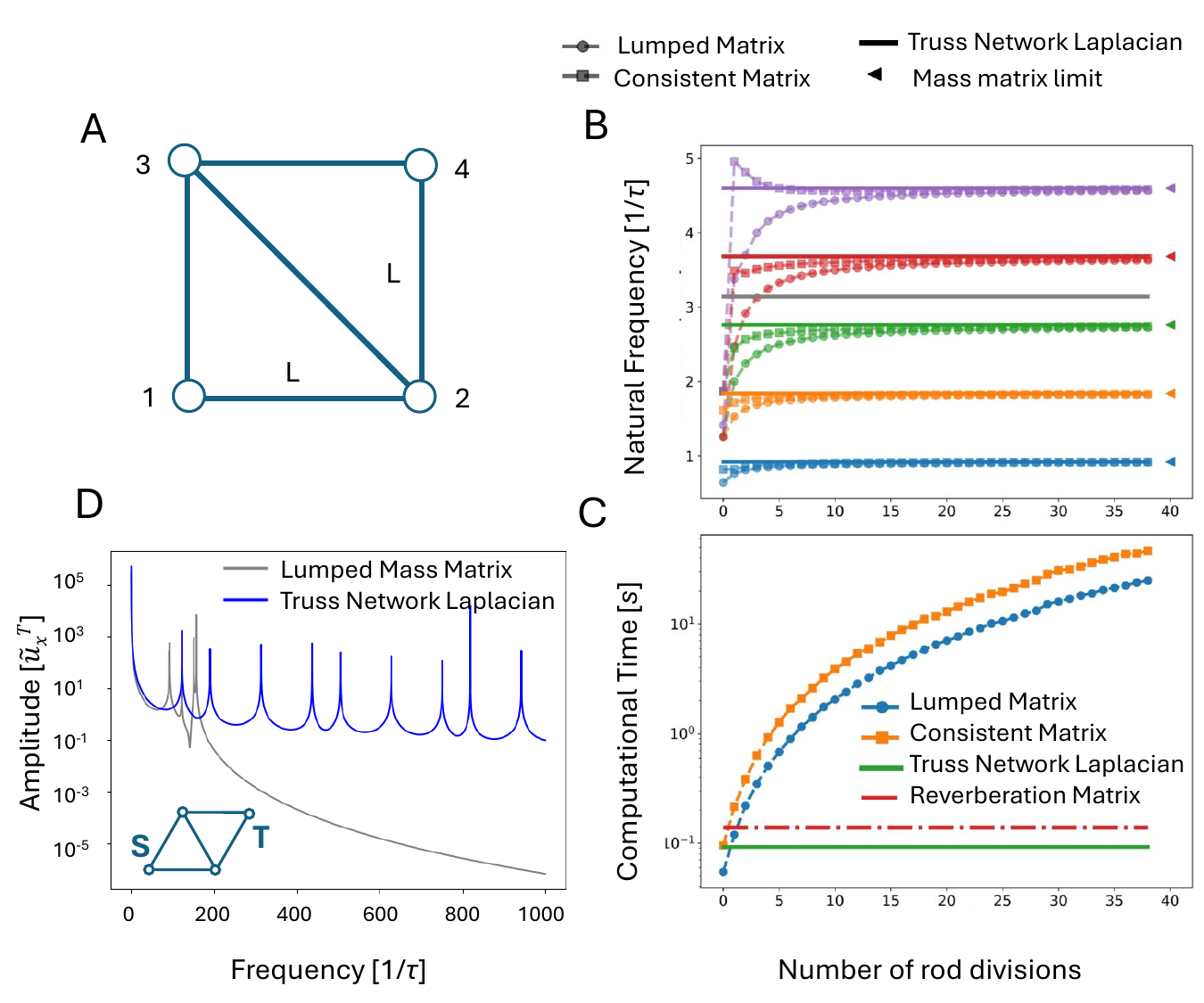}
    \caption{A. Schematic of a square truss structure with a crossbar. B. Comparison of natural frequencies computed using the truss network method and standard finite element methods for the square structure with a crossbar. The natural frequencies are plotted as a function of the number of subdivisions of the rods in the finite element methods. No subdivisions are necessary with the Network Laplacian method to achieve the theoretical limit. C. Computational time required for different methods. D. Amplitude of oscillation at the target joint (T) along the x-axis, $\tilde{u}_x^T$, in response to a periodic force applied at the source joint (S) along the (1,1) direction in a triangular structure with four joints (inset). Parameters: $\Lambda = 0.01$, $\tau = 0.01$.}
    \label{fig:truss-spring}
\end{figure*}
To show that our approach produces results that are equivalent to previously established methods, we first consider a two dimensional structure comprised of four joints arranged in a square, a total of 5 rods (square frame with a crossbar) as shown in Fig. \ref{fig:truss-spring} (A). We choose a global coordinate and label system such that the joint 1 is located at $(0,0)$, joint 2 at $(L,0)$, joint 3 at $(0,L)$, and joint 4 at $(L,L)$. Rods that form the square frame can be denoted by the connectivity labels 12, 13, 24, and 34. The fifth rod, labeled 23, creates a cross bar of length $\sqrt{2}L$. 

How this system responds to external forces can be determined from the use of scattering wavefronts (See supplementary information, Sec. S3). Here we focus on the natural frequencies given by our method and compare with finite element methods. We compare with two different finite element approaches. In finite element methods, one would define a mass matrix, $\mathbf{\overset\Leftrightarrow{M}}$, then determine the natural frequencies by observing where $\det(\mathbf{\overset\Leftrightarrow{K}}-\omega^{2}\mathbf{\overset\Leftrightarrow{M}})=0$. Eq. \ref{MCdef} gives one possible definition of $\mathbf{\overset\Leftrightarrow{M}}$, denoted as the consistent mass matrix, where we divide mass between rods in a way that is consistent with the network Laplacian approach. Another even simpler definition is known as the lumped mass matrix in which $N(x)=\Theta(1/2-x)$, with $\Theta(x)$ being the Heaviside step function, thus causing $\mathbf{\overset\Leftrightarrow{M}}$ to become diagonal with $\mathbf{\overset\Leftrightarrow{M}}_{\mu\mu}=(1/2)m_{\mu}\overset\Leftrightarrow{I}_{\mu}$. Here $m_{\mu}$ is the total mass of all rods connected to joint $\mu$. This lumped mass matrix simply compresses the entire system into the joints by splitting the mass of each element on the rod evenly between its connected joints. (Note that if the rods in the original network are not subdivided, the lumped mass matrix method is the same as the simple balls and springs method). The size of these matrices are determined by the number of subdivisions in the rods.

Either of these mass matrices can be used to calculate a set of natural frequencies, though these will generally not be equivalent due to the distinct way each handles the mass of the various system elements. The theory presented here is also distinct from this approach as the network Laplacian was developed without any approximation of the distribution of mass within the rods. We thus expect the natural frequencies calculated from the condition $\det(\mathbf{\overset\Leftrightarrow{D}})=0$ to represent a continuum limit of those derived from either mass matrix. To show this, we considered the effects of dividing each bar of the square structure shown in Fig. \ref{fig:truss-spring} (B) into a number of elements. Specifically, for $n$ divisions we introduce $n-1$ evenly spaced new joints in each rod to transform it into $n$ smaller identical rods. We then calculate the five lowest nonzero natural frequencies using all four methods (the lumped mass matrix, consistent mass matrix, network Laplacian, and reverberation matrix 
by steadily increasing $\omega$ and observing where the determinant of each matrix vanishes. The spectrum of the network Laplacian and reverberation matrix are completely equivalent, but the reverberation method ultimately requires more computational time.

Fig. \ref{fig:truss-spring} (B) shows the results of this process. The network Laplacian (and, equivalently, the reverberation matrix) has a unique representation for a given network and does not involve rod subdivisions. Hence the results from these methods are represented as horizontal lines. Note that the natural frequencies of both lumped mass matrix and consistent mass matrix methods steadily converge to those of the network Laplacian as the structure becomes more finely divided. This explicitly shows how our model acts as a continuum limit of the mass matrix method. The network Laplacian produces natural frequencies equivalent to the mass matrices in the limit of infinite subdivision of the system elements. In the particular case considered here in which all rods are made of identical material and have identical cross sections, the network Laplacian and reverberation matrix also find a natural mode at $\omega=\pi c/L$ that the mass matrices do not. This is due to this particular value of $\omega$ being a resonant frequency of each rod except for the cross bar, thus requiring the more careful handling explored in the supplementary information Sec. S1.

We also compare the overall computation times required to obtain the spectra plotted in Fig. \ref{fig:truss-spring} (C) as performed on an AMD Ryzen 5 1500x processor with Numpy's linear algebra package used to calculate the determinants. Fig. \ref{fig:truss-spring} (C) shows how these computation times increase for the mass matrix methods as the system becomes more finely divided. As mentioned previously, the natural frequencies and modes determined by the network Laplacian and reverberation matrix are invariant to such divisions. We see from this data that the network Laplacian is indeed more computationally efficient than the reverberation matrix method that involves larger matrices and that both are substantially more efficient than either mass matrix method when the system rods are subdivided into more than 3 divisions.\\

In many earlier studies on such network structures, a simple balls and springs method was used \cite{bergne2022scalable, sheinman2012nonlinear, rocks2017designing, huang2023jammed}. This method assumes that the mass of the rods is concentrated at the joints and the edges are replaced by a massless hookean spring. Note that this is same as the lumped mass matrix method with no sub-division of the edges. We computed the responses of the network with the truss method and balls and springs for a simple small triangular structure and also larger disordered networks that are of interest to soft matter research.

For balls and springs we assume the mass of the rods is equally shared between the two joints connecting them and is concentrated at the joints. The stiffness of the rod connecting joints $\mu$ and $\nu$ is given by $k_{\mu\nu}=\Lambda_{\mu\nu}/ \tau_{\mu \nu}$ while its mass is given by $m_{\mu\nu}=\Lambda_{\mu\nu} \tau_{\mu \nu}$. This for balls and springs matrix representation equivalent to a given truss network is:
\begin{equation}
\mathbf{\overset\Leftrightarrow{K}} - \mathbf{\omega^2 \overset\Leftrightarrow{M}} = \begin{cases}

\sum_{\gamma\in\mathcal{N}_{\mu}}k_{\mu\gamma}\hat{e}_{\mu\gamma}\hat{e}_{\mu\gamma}^{\text{T}}-\frac{m_{\mu\gamma}\omega^2}{2}\overset\Leftrightarrow{\mathrm{I}}_\mu& \mu=\nu \\
k_{\mu\nu}\hat{e}_{\mu\nu}\hat{e}_{\nu\mu}^{\text{T}} & \nu\in\mathcal{N}_{\mu} \\
\overset\Leftrightarrow{0} & \text{otherwise} \end{cases}.
\label{DBS}
\end{equation}

We first examined how our truss network-based method compares with balls and springs by examining computations for a triangular network with 4 joints with all rods having the same $\Lambda$ and $\tau$ for simplicity (Fig. \ref{fig:truss-spring} D (inset)). We applied a unit input force along (1, 1) to the source joint \textit{S} at (0,0) and measured the response at the diagonally opposite joint, \textit{T}. The response, defined as the amplitude of oscillation at the target joint, was analyzed as a function of frequency. At low frequencies, the truss and spring methods exhibited similar behaviors, consistent with our expectation that balls and springs are a good approximation for static and low-frequency cases. The peaks in the response correspond to resonant frequencies. Note that the location of the resonant frequencies differs between the two methods (which is in agreement with Fig. \ref{fig:truss-spring} (B)). More importantly, beyond a frequency of approximately 200, the spring network's response dropped rapidly, while the truss network continued to oscillate. This arises because a balls and spring network is a discrete representation and has only a finite number of normal modes. A spring network cannot respond accurately to excitations beyond the highest natural frequency (the Debye frequency). Truss networks treat rods as continuous objects, allowing each rod to support an infinite number of normal modes, unlike a single spring that supports only one mode. \\

\begin{figure}
    \centering
    \includegraphics[width=1\linewidth]{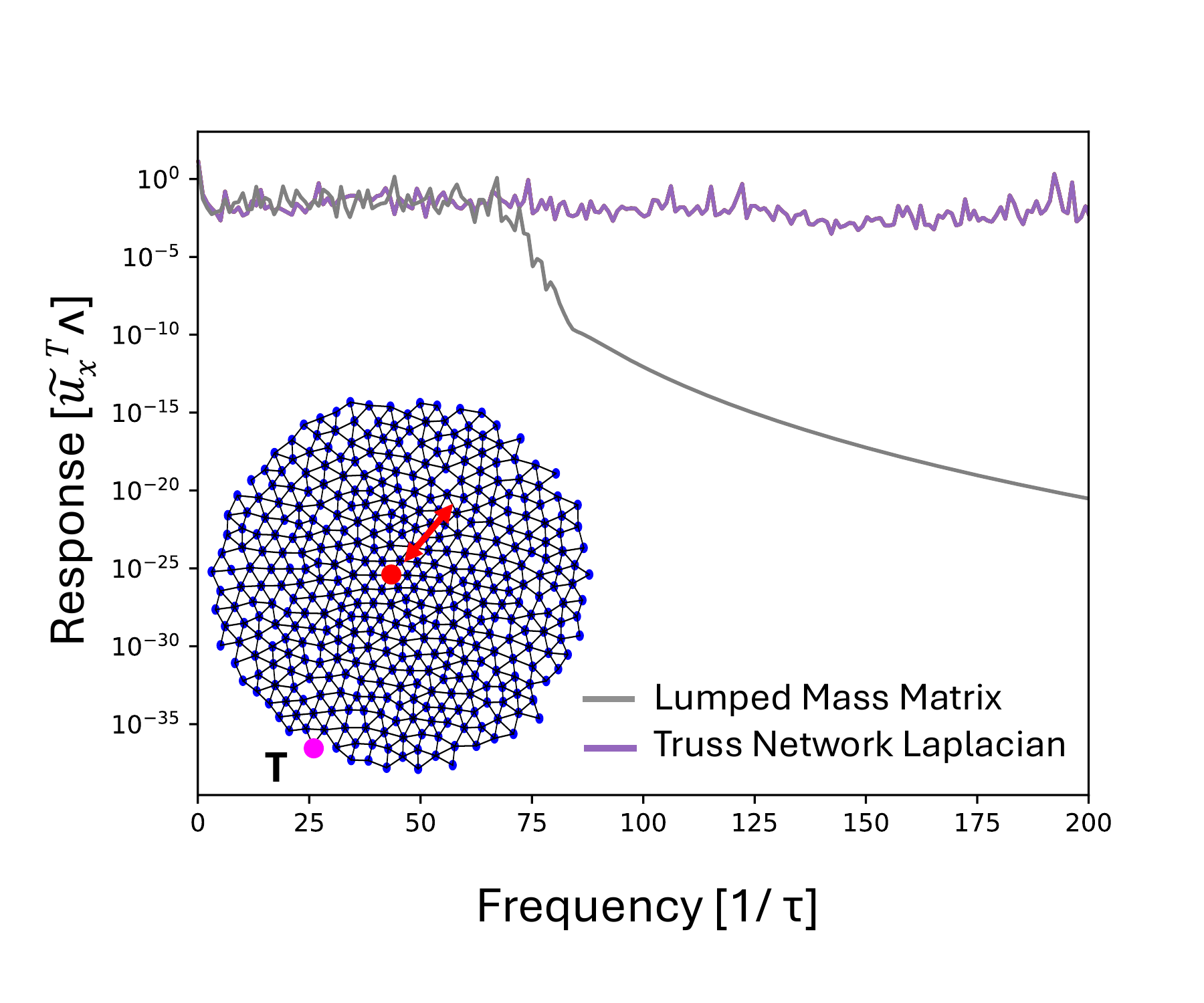}
    \caption{Comparison between the lumped mass matrix method (balls and springs) and the truss network method in a disordered network (inset). A periodic excitation is applied at the source joint (red) in the direction indicated by the red arrow. The oscillations at the target ($\tilde{u}_x^T$) are computed at the target joint (magenta).}
    \label{disordered_net}
\end{figure}

So far we illustrated our method using small regular networks. These are much simpler to compute, even allow for analytical calculations. Disordered networks are more relevant for metamaterial synthesis due to their high tunability and capacity to incorporate multiple functionalities \cite{rocks2017designing, rocks2019limits}. To evaluate our method in such contexts, we compared it with the lumped mass matrix (balls and springs) approach on a disordered network derived from jammed packing of spheres (Fig \ref{disordered_net}). For a disordered network with 400 joints (see inset of Fig \ref{disordered_net}) we compute the network response and report it as $\tilde{u}_x^T\Lambda$, where $\tilde{u}$ is the amplitude of oscillation at the target joint along x-axis for a unit force along (1,1) at the source joint. We found that beyond a certain frequency, the spring method fails to generate any response at the target joint. In contrast, the truss network method continues to yield consistent results (consistent with results discussed on a small network previously).
\subsection{Using our method to search for structures with enhanced vibrations at target}
\begin{figure*}
    \centering
    \includegraphics[width=\linewidth]{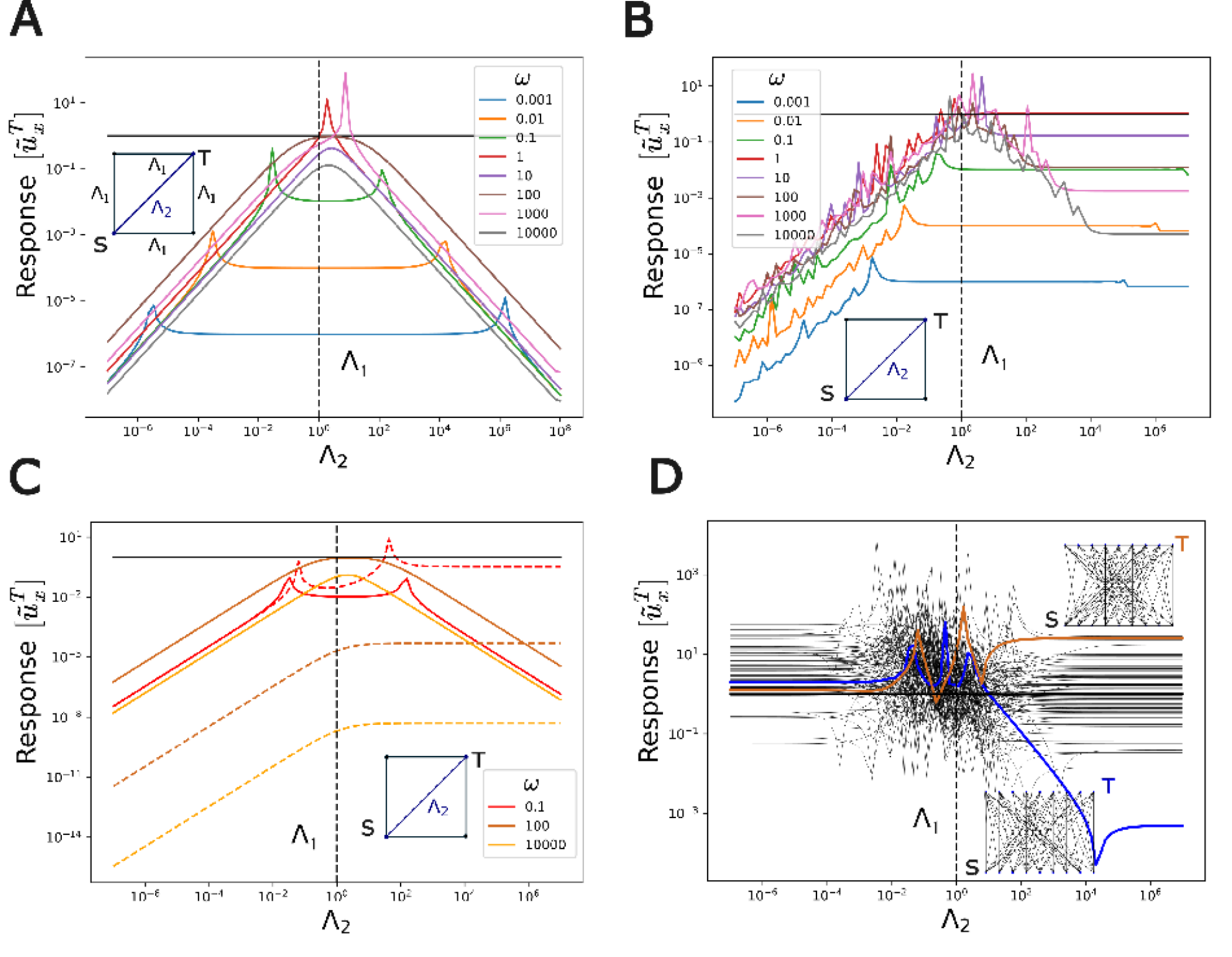}
    \caption{A. Tuning the response of a square network with a crossbar. Sinusoidal displacements of unit magnitude are given as input along x-axis at the source joint (S) and the amplitude of oscillation at the target joint (T) is computed. The other two diagonally opposite joints are given appropriate boundary conditions to prevent rigid body translation and rotation. The diagonal rod has an impedance $\Lambda_2$, rest of the rods have an impedance $\Lambda_1$. $\tau_1 = \tau_2 = 1$. B. Response of the target joint as a function of impedance keeping the total mass of the structure constant, by setting $\tau_2 = 1/ \Lambda_2$. , $\tau_1=1$. C. The amplitude of oscillation of the target joint (T) for different frequencies from the truss method (solid lines) compared with balls and spring (dashed lines). $\tau_1=\tau_2 =1$. D. Responses of random networks with slanted rod connections between two parallel rows of rods. Light grey lines represent responses from different network realizations. Source (S) is at (0,0), and target (T) is at (8,1.0). All slanted rods have impedance $\Lambda_2$, while horizontal rods have impedance $\Lambda_1$. Grey lines represent responses for different random realizations of the network. Blue and Chocolate colored lines indicate two selected realizations with markedly different responses; insets near the curves show corresponding network structures (shown in rectangular aspect ratio for legibility). The mass of the network was not held constant as $\Lambda_2$ was varied. Time constants and frequency are set as $\tau_1 = \tau_2 = 1$ and $\omega = 1$, respectively.}
    \label{fig:imedancematch}
\end{figure*}
The ability to control the oscillation of specific joints has broad applications, ranging from earthquake and blast protection systems to piezoelectric devices. Previous reports have discussed impedance-matching or mismatching protocols to achieve the desired control \cite{rathod2019review, rahimzadeh2015design, chen2014functionally, guo2021impedance, zhang2020engineering, verma2019impedance, ma2022enhancing}. While impedance matching in one-dimensional systems has been extensively studied \cite{kolsky1963stress}, two-dimensional and disordered networks present significantly more complex challenges \cite{blaurock1994impedance, zhu2022vibration}, requiring computational approaches. \\

Because our network Laplacian method is accurate and fast, it can be used to predict designs with novel mechanical prooperties. Here we sweep through design space of parameters and show structures with widely different oscillations of a target joint (T) for the same boundary conditions. We consider an example in which the target joint (T) oscillates in response to input oscillations at a source joint (S) within a square network with a crossbar, subject to boundary conditions that prevent rigid body translation and rotation. (Fixed boundary condition, $\tilde{u}_x = \tilde{u}_y =0$ for lower right joint, and roller boundary condition $\tilde{u}_y=0$ for upper left joint) (inset in Fig. \ref{fig:imedancematch} A). \\

We first analyzed the response of a square frame with a crossbar as a function of the impedance of crossbar ($\Lambda_2 = A\sqrt{E\rho}$) in  (Fig. \ref{fig:imedancematch} A). The four other rods have impedance, $\Lambda_1$. When $\Lambda_2$ is much smaller than $\Lambda_1$, the target joint's response is minimal, analogous to a thin rod ineffectively transmitting forces from the source joint. As $\Lambda_2$ increases, improved momentum transfer leads to an increased target response, reaching a maximum at a specific $\Lambda_2$ close to $\Lambda_2=\Lambda_1$. Note that this optimal impedance value is not simply equal to $\Lambda_1$ but is a complex function of the network architecture and the frequency of the force signal at the joint. (In contrast, the equivalent one-dimensional system has the oscillation amplitude maximized when $\Lambda_2 = \Lambda_1$, when $\tau_1=\tau_2=1$ for $
\omega =100$. See supplementary information, Sec. S 5.2, two rods end to end). We have derived an analytical relationship between response and parameters $\Lambda$ and oscillation frequency $\omega$ for a single rod, confirming the accuracy of our computational model (See supplementary information, Sec. S 5.1). In addition to maximized response when impedance is matched, there can also be some resonant modes, which lead to sharper peaks in this response.  

Further investigation involved modifying impedances $\Lambda_2$  (A$\sqrt{E\rho}$) while maintaining constant total mass by adjusting $\tau_2$ ( $L \sqrt{\frac{\rho}{E}}$) accordingly as $\tau_2$ = $1/\Lambda_2$ (Fig. \ref{fig:imedancematch} B). This approach revealed multiple amplitude maxima due to resonances, occurring when the $\cot(\omega \tau)$ term in the Laplacian matrix, $\mathbf{\overset\Leftrightarrow{D}}$ diverges. 

We then tested how this method compares with balls and springs methods (lumped mass matrix without any rod subdivisions) (Fig. \ref{fig:imedancematch} C). This analysis also confirms that the response of the network Laplacian method and the ball-and- spring method is comparable at lower frequencies but the ball-and-spring method is not adequate for higher frequencies as discussed in previous sections. It is interesting to note, at higher frequencies, the ball-and-spring method gives similar qualitative behavior as the network Laplacian method when $\Lambda_2$ is small relative to $\Lambda_1$, however at higher values of $\Lambda_2$, the ball-and-spring method gives a constant response while the network Laplacian method shows a decreased response. 

Next, we apply the method to disordered networks and give an example where our method shows that simple changes in network connectivity can lead to large changes in response (Fig. \ref{fig:imedancematch} D). We consider a disordered structure comprising 18 joints, equidistantly arranged in two parallel horizontal lines at $y = 0$ and $y = 1$. The network is constructed by randomly pairing joints at $y=0$ and $y=1$ and connecting them with slanted rods until the average coordination number reaches $z=4$, ensuring network rigidity. All horizontal rods have an impedance of $\Lambda_1 = 1$, while all slanted (and vertical) rods have a variable impedance $\Lambda_2$.
The boundary conditions are defined as follows: the source joint (lower leftmost joint (S) at (0,0)) is subjected to a constant unit displacement amplitude, $\tilde{u}_x^1=1$; the lower rightmost joint, at (8,0) is fixed along x and y axes, $\tilde{u}_x=\tilde{u}_y=0$; and the upper leftmost joint, at (0, 1) is fixed along the y-axis, $\tilde{u}_y=0$. These conditions prevent rigid body translations and rotations. We measure the oscillation amplitude (response) of the target joint (T) (upper rightmost, at (8,1.0)) as a function of $\Lambda_2$, with fixed parameters $\Lambda_1=1$, $\omega=1$, and $\tau_1 = \tau_2 =1$.
Fig. \ref{fig:imedancematch} (D) reveals that the target joint's response is highly sensitive to the specific network realization. The grey lines represent different random network configurations, all maintaining the same average coordination number but with unique pairings of joints at $y=0$ and $y=1$. The response curves exhibit three distinct regimes:
In the low $\Lambda_2$ regime ($\Lambda_2 \ll \Lambda_1$), the response is relatively constant as $\Lambda_2$ increases. This contrasts with the low and increasing response observed in simpler geometries (Fig. \ref{fig:imedancematch} A, C), attributable to the fact that in this system, even vertical rods have impedance of $\Lambda_2$, which is not the case for square with a crossbar system considered in Fig. \ref{fig:imedancematch}. Even a square network with a crossbar with a more number of $\Lambda_2$ elements shows similar behavior in this regime (See Fig. S6)
The intermediate $\Lambda_2$ regime ($\Lambda_2 \approx \Lambda_1$) shows many resonance peaks. Maximum variability between network realizations is observed here, highlighting the system's sensitivity to specific network realization.
In the high $\Lambda_2$ regime ($\Lambda_2 \gg \Lambda_1$), the response reaches a plateau (but distinct for each network realization). We identified two network realizations with markedly different responses in the high $\Lambda_2$ regime ($\Lambda_2 \gg \Lambda_1$). The corresponding network structure is included as insets near these curves. This shows that subtle changes in network structures can lead to several orders of magnitude differences in the responses. However, it is to be noted that whole range of $\Lambda_2$ values considered in this study may not be accessible by changing areas of cross-section of rods. Larger areas could lead to steric hindrance with neighbours. Nonetheless, the strategy of changing impedances with changes in area of cross section offers a valuable and practical design approach for 3D-printed microarchitectured materials to tune their response.\\

\section{Conclusions}
In this work, we developed a computationally tractable approach for computing the dynamic responses of complex truss networks. We derived a network Laplacian formalism that relates displacements (or velocities) of joints to applied forces, taking into account the mass distribution along the rods. We demonstrated that our method accurately computes natural frequencies of truss networks through comparison with more detailed finite element methods on simple structures. This method is at least an order of magnitude faster than mass matrix techniques (finite element) for a reasonable level of accuracy. We then compared with balls and springs model and showed that while balls and springs model gives comparable results at low frequencies, it breaks down at high frequencies (beyond the highest normal mode or the Debye frequency of spring networks). \\

Our method can contribute to research efforts in uncovering the effects of network topologies, impedance distributions, and other network properties on the mechanical responses of disordered soft matter systems. For example, our specific observation of high variability between random realizations in Fig.~\ref{fig:imedancematch}(D) illustrates the potential for fine-tuning mechanical properties through structural design using our method. Moreover, the contrast with simpler geometries emphasizes the importance of considering larger networks in predicting and designing mechanical responses. Furthermore, the observed resonance behavior in the intermediate $\Lambda_2$ regime hints at the possibility of designing switchable materials, where small changes in the network properties could lead to significant shifts in mechanical response. This could have applications in areas such as vibration damping, acoustic metamaterials, or responsive soft robotic systems. With suitable optimization techniques, this approach could aid in designing materials with tailored mechanical properties through the large design space provided by disordered networks.\\

We believe this approach could be a valuable addition to computational toolkits for designing micro-architectured 3D printable materials. This method, coupled with suitable optimization methods, can help in efficiently designing networks with nontrivial mechanical properties.  One can also apply the method to understand more complex disordered and hierarchical biological structures, such as bones. Finite element methods are not appropriate for hierarchical cases as they require refinement of mesh sizes when modeling different levels of arrangement. The truss network Laplacian method is specifically useful for hierarchical truss structures as it relies only on joint positions and connectivities and doesn't involve the discretization of rods.\\

While our method is efficient and runs fast for a wide range of cases it does suffer from several drawbacks that in some contexts limit its applicability.  Currently, our approach is appropriate for designing small-amplitude dynamics of stretch-dominated structures. We ignore transverse forces, \textit{i.e.}, rods cannot transmit forces in directions perpendicular to their long axis. We also ignore the bending of the rods. This can be crucial if the structure is floppy or if the applied forces are large. Thus, this method, while universally applicable to any network of rods with pin joints, becomes particularly powerful and accurate for 3D printed materials when coordination numbers are at or above the critical value dictated by Maxwell's criterion for rigidity. Finally, we considered perfectly elastic rods that don't dissipate any mechanical energy. This is not true for real materials, especially polymeric materials used for 3D printing which dissipate mechanical energy over time. \\ 

Future directions include further improving the model by incorporating rod bending, transverse forces, and viscoelasticity. These enhancements would make the model more suitable for designing networks with complex mechanical properties for real-world applications. In particular, this approach could be used to design metamaterials with controlled mechanical energy - an area of interest for soft robotics. Additionally, the method can help study the stress response of complex hierarchical structures, such as bone tissue, and how these tissues remodel under stress. This may offer insights into adaptation principles under external periodic forces and aligns with current efforts in soft matter research to explore physical learning in biological systems. Beyond understanding biological structures, this approach can also be applied to the development of artificial adaptive mechanical metamaterials that dynamically modify their properties in response to external conditions. \\



\section*{Code availability}
Code from this method is freely available at \url{https://github.com/nsarpangala/truss_network_spectral_node}
\section*{Author Contributions}
SF and NS contributed equally to the paper. SF contributed to model development, running numerical computations, analysing data, and writing the paper. NS worked on developing the code, running numerical computations, analysing data and writing paper. PKP and EK conceived the study, developed the formalism, and wrote the paper. 

\section*{Conflicts of interest}
There are no conflicts to declare.

\section*{Acknowledgements}
This research was funded by the ARO MURI grant 10085212, the  University  of  Pennsylvania  Materials Research Science and Engineering Center (MRSEC) through Grant No. DMR-1720530 and DMR-2309043 and the Simons Foundation through Grant No. 568888. PKP acknowledges support for this work through a seed grant from Penn's Materials Science and Engineering Center (MRSEC) grant DMR-1720530.



\balance


\bibliography{refs} 

\providecommand*{\mcitethebibliography}{\thebibliography}
\csname @ifundefined\endcsname{endmcitethebibliography}
{\let\endmcitethebibliography\endthebibliography}{}
\begin{mcitethebibliography}{38}
\providecommand*{\natexlab}[1]{#1}
\providecommand*{\mciteSetBstSublistMode}[1]{}
\providecommand*{\mciteSetBstMaxWidthForm}[2]{}
\providecommand*{\mciteBstWouldAddEndPuncttrue}
  {\def\EndOfBibitem{\unskip.}}
\providecommand*{\mciteBstWouldAddEndPunctfalse}
  {\let\EndOfBibitem\relax}
\providecommand*{\mciteSetBstMidEndSepPunct}[3]{}
\providecommand*{\mciteSetBstSublistLabelBeginEnd}[3]{}
\providecommand*{\EndOfBibitem}{}
\mciteSetBstSublistMode{f}
\mciteSetBstMaxWidthForm{subitem}
{(\emph{\alph{mcitesubitemcount}})}
\mciteSetBstSublistLabelBeginEnd{\mcitemaxwidthsubitemform\space}
{\relax}{\relax}

\bibitem[Askari \emph{et~al.}(2020)Askari, Hutchins, Thomas, Astolfi, Watson, Abdi, Ricci, Laureti, Nie, Freear,\emph{et~al.}]{askari2020additive}
M.~Askari, D.~A. Hutchins, P.~J. Thomas, L.~Astolfi, R.~L. Watson, M.~Abdi, M.~Ricci, S.~Laureti, L.~Nie, S.~Freear \emph{et~al.}, \emph{Additive Manufacturing}, 2020, \textbf{36}, 101562\relax
\mciteBstWouldAddEndPuncttrue
\mciteSetBstMidEndSepPunct{\mcitedefaultmidpunct}
{\mcitedefaultendpunct}{\mcitedefaultseppunct}\relax
\EndOfBibitem
\bibitem[Torres \emph{et~al.}(2019)Torres, Trikanad, Aubin, Lambers, Luna, Rimnac, Zavattieri, and Hernandez]{torres2019bone}
A.~M. Torres, A.~A. Trikanad, C.~A. Aubin, F.~M. Lambers, M.~Luna, C.~M. Rimnac, P.~Zavattieri and C.~J. Hernandez, \emph{Proceedings of the National Academy of Sciences}, 2019, \textbf{116}, 24457--24462\relax
\mciteBstWouldAddEndPuncttrue
\mciteSetBstMidEndSepPunct{\mcitedefaultmidpunct}
{\mcitedefaultendpunct}{\mcitedefaultseppunct}\relax
\EndOfBibitem
\bibitem[Van~Tol \emph{et~al.}(2020)Van~Tol, Schemenz, Wagermaier, Roschger, Razi, Vitienes, Fratzl, Willie, and Weinkamer]{van2020mechanoresponse}
A.~F. Van~Tol, V.~Schemenz, W.~Wagermaier, A.~Roschger, H.~Razi, I.~Vitienes, P.~Fratzl, B.~M. Willie and R.~Weinkamer, \emph{Proceedings of the National Academy of Sciences}, 2020, \textbf{117}, 32251--32259\relax
\mciteBstWouldAddEndPuncttrue
\mciteSetBstMidEndSepPunct{\mcitedefaultmidpunct}
{\mcitedefaultendpunct}{\mcitedefaultseppunct}\relax
\EndOfBibitem
\bibitem[Jung \emph{et~al.}(2018)Jung, Pissarenko, Yaraghi, Naleway, Kisailus, Meyers, and McKittrick]{jung2018comparative}
J.-Y. Jung, A.~Pissarenko, N.~A. Yaraghi, S.~E. Naleway, D.~Kisailus, M.~A. Meyers and J.~McKittrick, \emph{Journal of the Mechanical Behavior of Biomedical Materials}, 2018, \textbf{84}, 273--280\relax
\mciteBstWouldAddEndPuncttrue
\mciteSetBstMidEndSepPunct{\mcitedefaultmidpunct}
{\mcitedefaultendpunct}{\mcitedefaultseppunct}\relax
\EndOfBibitem
\bibitem[Mueller \emph{et~al.}(2019)Mueller, Matlack, Shea, and Daraio]{mueller2019energy}
J.~Mueller, K.~H. Matlack, K.~Shea and C.~Daraio, \emph{Advanced Theory and Simulations}, 2019, \textbf{2}, 1900081\relax
\mciteBstWouldAddEndPuncttrue
\mciteSetBstMidEndSepPunct{\mcitedefaultmidpunct}
{\mcitedefaultendpunct}{\mcitedefaultseppunct}\relax
\EndOfBibitem
\bibitem[Glaesener \emph{et~al.}(2021)Glaesener, Bastek, Gonon, Kannan, Telgen, Sp{\"o}ttling, Steiner, and Kochmann]{glaesener2021viscoelastic}
R.~N. Glaesener, J.-H. Bastek, F.~Gonon, V.~Kannan, B.~Telgen, B.~Sp{\"o}ttling, S.~Steiner and D.~M. Kochmann, \emph{Journal of the Mechanics and Physics of Solids}, 2021, \textbf{156}, 104569\relax
\mciteBstWouldAddEndPuncttrue
\mciteSetBstMidEndSepPunct{\mcitedefaultmidpunct}
{\mcitedefaultendpunct}{\mcitedefaultseppunct}\relax
\EndOfBibitem
\bibitem[Wang \emph{et~al.}(2015)Wang, Casadei, Kang, and Bertoldi]{wang2015locally}
P.~Wang, F.~Casadei, S.~H. Kang and K.~Bertoldi, \emph{Physical Review B}, 2015, \textbf{91}, 020103\relax
\mciteBstWouldAddEndPuncttrue
\mciteSetBstMidEndSepPunct{\mcitedefaultmidpunct}
{\mcitedefaultendpunct}{\mcitedefaultseppunct}\relax
\EndOfBibitem
\bibitem[Krushynska \emph{et~al.}(2018)Krushynska, Amendola, Bosia, Daraio, Pugno, and Fraternali]{krushynska2018accordion}
A.~O. Krushynska, A.~Amendola, F.~Bosia, C.~Daraio, N.~M. Pugno and F.~Fraternali, \emph{New Journal of Physics}, 2018, \textbf{20}, 073051\relax
\mciteBstWouldAddEndPuncttrue
\mciteSetBstMidEndSepPunct{\mcitedefaultmidpunct}
{\mcitedefaultendpunct}{\mcitedefaultseppunct}\relax
\EndOfBibitem
\bibitem[Zhang \emph{et~al.}(2015)Zhang, Koo, Liu, Zou, Chattopadhyay, and Dai]{zhang2015novel}
J.~Zhang, B.~Koo, Y.~Liu, J.~Zou, A.~Chattopadhyay and L.~Dai, \emph{Smart Materials and Structures}, 2015, \textbf{24}, 085022\relax
\mciteBstWouldAddEndPuncttrue
\mciteSetBstMidEndSepPunct{\mcitedefaultmidpunct}
{\mcitedefaultendpunct}{\mcitedefaultseppunct}\relax
\EndOfBibitem
\bibitem[Lu \emph{et~al.}(2022)Lu, Liu, and Hu]{lu2022double}
Y.~Lu, X.-Y. Liu and G.-H. Hu, \emph{Macromolecules}, 2022, \textbf{55}, 4548--4556\relax
\mciteBstWouldAddEndPuncttrue
\mciteSetBstMidEndSepPunct{\mcitedefaultmidpunct}
{\mcitedefaultendpunct}{\mcitedefaultseppunct}\relax
\EndOfBibitem
\bibitem[Branch \emph{et~al.}(2017)Branch, Ionita, Clements, Montgomery, Jensen, Patterson, Schmalzer, Mueller, and Dattelbaum]{branch2017controlling}
B.~Branch, A.~Ionita, B.~E. Clements, D.~S. Montgomery, B.~J. Jensen, B.~Patterson, A.~Schmalzer, A.~Mueller and D.~M. Dattelbaum, \emph{Journal of Applied Physics}, 2017, \textbf{121}, 135102\relax
\mciteBstWouldAddEndPuncttrue
\mciteSetBstMidEndSepPunct{\mcitedefaultmidpunct}
{\mcitedefaultendpunct}{\mcitedefaultseppunct}\relax
\EndOfBibitem
\bibitem[Cook \emph{et~al.}(2007)Cook\emph{et~al.}]{cook2007concepts}
R.~D. Cook \emph{et~al.}, \emph{Concepts and applications of finite element analysis}, John wiley \& sons, 2007\relax
\mciteBstWouldAddEndPuncttrue
\mciteSetBstMidEndSepPunct{\mcitedefaultmidpunct}
{\mcitedefaultendpunct}{\mcitedefaultseppunct}\relax
\EndOfBibitem
\bibitem[Howard and Pao(1998)]{howard1998analysis}
S.~M. Howard and Y.-H. Pao, \emph{Journal of Engineering Mechanics}, 1998, \textbf{124}, 884--891\relax
\mciteBstWouldAddEndPuncttrue
\mciteSetBstMidEndSepPunct{\mcitedefaultmidpunct}
{\mcitedefaultendpunct}{\mcitedefaultseppunct}\relax
\EndOfBibitem
\bibitem[Pao \emph{et~al.}(1999)Pao, Keh, and Howard]{pao1999dynamic}
Y.-H. Pao, D.-C. Keh and S.~M. Howard, \emph{AIAA journal}, 1999, \textbf{37}, 594--603\relax
\mciteBstWouldAddEndPuncttrue
\mciteSetBstMidEndSepPunct{\mcitedefaultmidpunct}
{\mcitedefaultendpunct}{\mcitedefaultseppunct}\relax
\EndOfBibitem
\bibitem[P{\"o}lz \emph{et~al.}(2019)P{\"o}lz, Gfrerer, and Schanz]{polz2019wave}
D.~P{\"o}lz, M.~H. Gfrerer and M.~Schanz, \emph{Wave Motion}, 2019, \textbf{87}, 37--57\relax
\mciteBstWouldAddEndPuncttrue
\mciteSetBstMidEndSepPunct{\mcitedefaultmidpunct}
{\mcitedefaultendpunct}{\mcitedefaultseppunct}\relax
\EndOfBibitem
\bibitem[Messner \emph{et~al.}(2015)Messner, Barham, Kumar, and Barton]{messner2015wave}
M.~C. Messner, M.~I. Barham, M.~Kumar and N.~R. Barton, \emph{International Journal of Solids and Structures}, 2015, \textbf{73}, 55--66\relax
\mciteBstWouldAddEndPuncttrue
\mciteSetBstMidEndSepPunct{\mcitedefaultmidpunct}
{\mcitedefaultendpunct}{\mcitedefaultseppunct}\relax
\EndOfBibitem
\bibitem[Trainiti \emph{et~al.}(2016)Trainiti, Rimoli, and Ruzzene]{trainiti2016wave}
G.~Trainiti, J.~J. Rimoli and M.~Ruzzene, \emph{International Journal of Solids and Structures}, 2016, \textbf{97}, 431--444\relax
\mciteBstWouldAddEndPuncttrue
\mciteSetBstMidEndSepPunct{\mcitedefaultmidpunct}
{\mcitedefaultendpunct}{\mcitedefaultseppunct}\relax
\EndOfBibitem
\bibitem[Fancher and Katifori(2022)]{fancher2022mechanical}
S.~Fancher and E.~Katifori, \emph{Physical Review Fluids}, 2022, \textbf{7}, 013101\relax
\mciteBstWouldAddEndPuncttrue
\mciteSetBstMidEndSepPunct{\mcitedefaultmidpunct}
{\mcitedefaultendpunct}{\mcitedefaultseppunct}\relax
\EndOfBibitem
\bibitem[JF(1997)]{jf1997wave}
D.~JF, \emph{Wave propagation in structures: spectral analysis using fast discrete Fourier transforms}, 1997\relax
\mciteBstWouldAddEndPuncttrue
\mciteSetBstMidEndSepPunct{\mcitedefaultmidpunct}
{\mcitedefaultendpunct}{\mcitedefaultseppunct}\relax
\EndOfBibitem
\bibitem[Chakraborty and Gopalakrishnan(2003)]{chakraborty2003spectrally}
A.~Chakraborty and S.~Gopalakrishnan, \emph{International Journal of Solids and Structures}, 2003, \textbf{40}, 2421--2448\relax
\mciteBstWouldAddEndPuncttrue
\mciteSetBstMidEndSepPunct{\mcitedefaultmidpunct}
{\mcitedefaultendpunct}{\mcitedefaultseppunct}\relax
\EndOfBibitem
\bibitem[Lee(2009)]{lee2009spectral}
U.~Lee, \emph{Spectral element method in structural dynamics}, John Wiley \& Sons, 2009\relax
\mciteBstWouldAddEndPuncttrue
\mciteSetBstMidEndSepPunct{\mcitedefaultmidpunct}
{\mcitedefaultendpunct}{\mcitedefaultseppunct}\relax
\EndOfBibitem
\bibitem[An \emph{et~al.}(2020)An, Fan, and Zhang]{an2020elastic}
X.~An, H.~Fan and C.~Zhang, \emph{Journal of sound and vibration}, 2020, \textbf{475}, 115292\relax
\mciteBstWouldAddEndPuncttrue
\mciteSetBstMidEndSepPunct{\mcitedefaultmidpunct}
{\mcitedefaultendpunct}{\mcitedefaultseppunct}\relax
\EndOfBibitem
\bibitem[Zuo \emph{et~al.}(2016)Zuo, Li, and Zhang]{zuo2016numerical}
S.-L. Zuo, F.-M. Li and C.~Zhang, \emph{Acta Mechanica}, 2016, \textbf{227}, 1653--1669\relax
\mciteBstWouldAddEndPuncttrue
\mciteSetBstMidEndSepPunct{\mcitedefaultmidpunct}
{\mcitedefaultendpunct}{\mcitedefaultseppunct}\relax
\EndOfBibitem
\bibitem[Bergne \emph{et~al.}(2022)Bergne, Baardink, Loukaides, and Souslov]{bergne2022scalable}
A.~Bergne, G.~Baardink, E.~G. Loukaides and A.~Souslov, \emph{Extreme Mechanics Letters}, 2022, \textbf{57}, 101911\relax
\mciteBstWouldAddEndPuncttrue
\mciteSetBstMidEndSepPunct{\mcitedefaultmidpunct}
{\mcitedefaultendpunct}{\mcitedefaultseppunct}\relax
\EndOfBibitem
\bibitem[Sheinman \emph{et~al.}(2012)Sheinman, Broedersz, and MacKintosh]{sheinman2012nonlinear}
M.~Sheinman, C.~Broedersz and F.~MacKintosh, \emph{Physical Review E—Statistical, Nonlinear, and Soft Matter Physics}, 2012, \textbf{85}, 021801\relax
\mciteBstWouldAddEndPuncttrue
\mciteSetBstMidEndSepPunct{\mcitedefaultmidpunct}
{\mcitedefaultendpunct}{\mcitedefaultseppunct}\relax
\EndOfBibitem
\bibitem[Rocks \emph{et~al.}(2017)Rocks, Pashine, Bischofberger, Goodrich, Liu, and Nagel]{rocks2017designing}
J.~W. Rocks, N.~Pashine, I.~Bischofberger, C.~P. Goodrich, A.~J. Liu and S.~R. Nagel, \emph{Proceedings of the National Academy of Sciences}, 2017, \textbf{114}, 2520--2525\relax
\mciteBstWouldAddEndPuncttrue
\mciteSetBstMidEndSepPunct{\mcitedefaultmidpunct}
{\mcitedefaultendpunct}{\mcitedefaultseppunct}\relax
\EndOfBibitem
\bibitem[Huang \emph{et~al.}(2023)Huang, Zhang, Xu, Zhang, Tong, and Xu]{huang2023jammed}
J.~Huang, J.~Zhang, D.~Xu, S.~Zhang, H.~Tong and N.~Xu, \emph{Current Opinion in Solid State and Materials Science}, 2023, \textbf{27}, 101053\relax
\mciteBstWouldAddEndPuncttrue
\mciteSetBstMidEndSepPunct{\mcitedefaultmidpunct}
{\mcitedefaultendpunct}{\mcitedefaultseppunct}\relax
\EndOfBibitem
\bibitem[Rocks \emph{et~al.}(2019)Rocks, Ronellenfitsch, Liu, Nagel, and Katifori]{rocks2019limits}
J.~W. Rocks, H.~Ronellenfitsch, A.~J. Liu, S.~R. Nagel and E.~Katifori, \emph{Proceedings of the National Academy of Sciences}, 2019, \textbf{116}, 2506--2511\relax
\mciteBstWouldAddEndPuncttrue
\mciteSetBstMidEndSepPunct{\mcitedefaultmidpunct}
{\mcitedefaultendpunct}{\mcitedefaultseppunct}\relax
\EndOfBibitem
\bibitem[Rathod(2019)]{rathod2019review}
V.~T. Rathod, \emph{Electronics}, 2019, \textbf{8}, 169\relax
\mciteBstWouldAddEndPuncttrue
\mciteSetBstMidEndSepPunct{\mcitedefaultmidpunct}
{\mcitedefaultendpunct}{\mcitedefaultseppunct}\relax
\EndOfBibitem
\bibitem[Rahimzadeh \emph{et~al.}(2015)Rahimzadeh, Arruda, and Thouless]{rahimzadeh2015design}
T.~Rahimzadeh, E.~M. Arruda and M.~Thouless, \emph{Journal of the Mechanics and Physics of Solids}, 2015, \textbf{85}, 98--111\relax
\mciteBstWouldAddEndPuncttrue
\mciteSetBstMidEndSepPunct{\mcitedefaultmidpunct}
{\mcitedefaultendpunct}{\mcitedefaultseppunct}\relax
\EndOfBibitem
\bibitem[Chen \emph{et~al.}(2014)Chen, Zhang, Hao, Lin, and Fu]{chen2014functionally}
S.~Chen, Y.~Zhang, C.~Hao, S.~Lin and Z.~Fu, \emph{Physics Letters A}, 2014, \textbf{378}, 77--81\relax
\mciteBstWouldAddEndPuncttrue
\mciteSetBstMidEndSepPunct{\mcitedefaultmidpunct}
{\mcitedefaultendpunct}{\mcitedefaultseppunct}\relax
\EndOfBibitem
\bibitem[Guo \emph{et~al.}(2021)Guo, Zheng, and Palffy-Muhoray]{guo2021impedance}
T.~Guo, X.~Zheng and P.~Palffy-Muhoray, \emph{Soft Matter}, 2021, \textbf{17}, 4191--4194\relax
\mciteBstWouldAddEndPuncttrue
\mciteSetBstMidEndSepPunct{\mcitedefaultmidpunct}
{\mcitedefaultendpunct}{\mcitedefaultseppunct}\relax
\EndOfBibitem
\bibitem[Zhang \emph{et~al.}(2020)Zhang, Qu, and Wang]{zhang2020engineering}
X.~Zhang, Z.~Qu and H.~Wang, \emph{Iscience}, 2020, \textbf{23}, xx--yy\relax
\mciteBstWouldAddEndPuncttrue
\mciteSetBstMidEndSepPunct{\mcitedefaultmidpunct}
{\mcitedefaultendpunct}{\mcitedefaultseppunct}\relax
\EndOfBibitem
\bibitem[Verma \emph{et~al.}(2019)Verma, Sivaselvan, and Rajasankar]{verma2019impedance}
M.~Verma, M.~Sivaselvan and J.~Rajasankar, \emph{Structural Control and Health Monitoring}, 2019, \textbf{26}, e2402\relax
\mciteBstWouldAddEndPuncttrue
\mciteSetBstMidEndSepPunct{\mcitedefaultmidpunct}
{\mcitedefaultendpunct}{\mcitedefaultseppunct}\relax
\EndOfBibitem
\bibitem[Ma \emph{et~al.}(2022)Ma, Wang, Du, Zhu, and Wu]{ma2022enhancing}
F.~Ma, C.~Wang, Y.~Du, Z.~Zhu and J.~H. Wu, \emph{Materials Horizons}, 2022, \textbf{9}, 653--662\relax
\mciteBstWouldAddEndPuncttrue
\mciteSetBstMidEndSepPunct{\mcitedefaultmidpunct}
{\mcitedefaultendpunct}{\mcitedefaultseppunct}\relax
\EndOfBibitem
\bibitem[Kolsky(1963)]{kolsky1963stress}
H.~Kolsky, \emph{Stress waves in solids}, Courier Corporation, 1963, vol. 1098\relax
\mciteBstWouldAddEndPuncttrue
\mciteSetBstMidEndSepPunct{\mcitedefaultmidpunct}
{\mcitedefaultendpunct}{\mcitedefaultseppunct}\relax
\EndOfBibitem
\bibitem[Blaurock(1994)]{blaurock1994impedance}
C.~A. Blaurock, \emph{Ph.D. thesis}, Massachusetts Institute of Technology, 1994\relax
\mciteBstWouldAddEndPuncttrue
\mciteSetBstMidEndSepPunct{\mcitedefaultmidpunct}
{\mcitedefaultendpunct}{\mcitedefaultseppunct}\relax
\EndOfBibitem
\bibitem[Zhu \emph{et~al.}(2022)Zhu, Wu, and Sun]{zhu2022vibration}
H.-Z. Zhu, J.-H. Wu and Y.-D. Sun, \emph{Applied Sciences}, 2022, \textbf{12}, 8863\relax
\mciteBstWouldAddEndPuncttrue
\mciteSetBstMidEndSepPunct{\mcitedefaultmidpunct}
{\mcitedefaultendpunct}{\mcitedefaultseppunct}\relax
\EndOfBibitem
\end{mcitethebibliography}
\bibliographystyle{rsc} 

\end{document}